\documentclass[11pt]{article}
\usepackage{amsthm}
\usepackage{amssymb}
\usepackage{amsbsy}
\usepackage{amsfonts}
\usepackage{amsmath}
\usepackage{amsopn}
\usepackage[active]{srcltx}

\newtheorem{theorem}{Theorem}
\newtheorem{proposition}[theorem]{Proposition}
\newtheorem{lemma}[theorem]{Lemma}

\DeclareMathOperator{\nr}{\mathbb{N}}
\DeclareMathOperator{\zr}{\mathbb{Z}}
\newcommand{\Z}{\mathbb{Z}}
\DeclareMathOperator{\re}{\mathbb{R}}
\newcommand{\R}{\mathbb{R}}
\DeclareMathOperator{\p}{\mathbb{P}}
\DeclareMathOperator{\e}{\mathbb{E}}

\DeclareMathOperator{\F}{\mathcal F}

\newcommand{\as}{a.s.}
\newcommand{\ie}{i.e.}
\newcommand{\eg}{e.g.}
\newcommand{\td}{\tilde}

\begin{document}
\author{A.~M.~G.~Cox\thanks{e-mail:
        \texttt{A.M.G.Cox@bath.ac.uk};
        web: \texttt{www.maths.bath.ac.uk/$\sim$mapamgc/}
\newline Author is grateful for financial support from the Nuffield Foundation}\\
        Dept.\ of Mathematical Sciences\\
        University of Bath\\
        Bath BA2 7AY, UK
%        Department of Mathematics\\
%        University of York\\ Heslington\\
%        York Y010 5DD, U.~K.
 \and Jan Ob\l \'oj\thanks{e-mail:
        \texttt{obloj@ccr.jussieu.fr}; web:
        \texttt{www.proba.jussieu.fr/$\sim$obloj/}}\\
        Laboratoire de Probabilit\'es\\
        Universit\'e Paris 6\\
        4 pl. Jussieu -- Boite 188\\75252 Paris Cedex 05, France}

%\author{A.M.G.\ Cox, Jan Ob\l \'oj}
\date{}
\title{Classes of Skorokhod Embeddings for the Simple Symmetric Random Walk}
\maketitle
\begin{abstract}
The Skorokhod Embedding problem is well understood when the underlying process is a Brownian motion. We examine the problem when the underlying is the simple symmetric random walk and when no external randomisation is allowed. We prove that any measure on $\zr$ can be embedded by means of a minimal stopping time. However, in sharp contrast to the Brownian setting, we show that the set of measures which can be embedded in a uniformly integrable way is strictly smaller then the set of centered probability measures: specifically it is a fractal set which we characterise as an iterated function system. Finally, we define the natural extension of several known constructions from the Brownian setting and show that these constructions require us to further restrict the sets of target laws.

\end{abstract}
\noindent \emph{2000 Mathematics Subject Classification}: \smallskip

\section{Introduction}

The Skorokhod embedding problem was first posed (and solved) in Skorokhod \cite{MR32:3082b}. Since then, the problem has been an active field of research and has found numerous solutions. We refer the reader to Ob\l\'oj \cite{genealogia} for a comprehensive survey paper.

Simply stated the problem is the following: given a probability measure $\mu$ and a stochastic process $(X_t)_{t \ge 0}$ find a stopping time $T$ such that $X_T\sim \mu$. The most commonly considered case is when $(X_t)$ is a 1-dimensional Brownian motion. However, in this context we have a trivial solution (usually attributed to Doob): for any probability measure $\mu$, define the distribution function $F_\mu(x)=\mu((-\infty,x])$ with $F_\mu^{-1}$ its right-continuous inverse, and let $\Phi$ denote the distribution function of a standard Normal variable. Then the stopping time $T=\inf\{t\geq 2: B_t=F_\mu^{-1}(\Phi(B_1))\}$ embeds $\mu$ in Brownian motion (\ie{} $B_T\sim \mu$). Thus it is clear that interest lies in the properties of the stopping time $T$.

In the example above we always have $\e T=\infty$. Skorokhod \cite{MR32:3082b} imposed $\e T<\infty$ which then implies that $\mu$ is centered with finite second moment and that the process $(B_{t\land T}:t\ge 0)$ is a uniformly integrable martingale. Numerous authors (e.g.\ Root \cite{MR38:6670}, Az\'ema and Yor \cite{MR82c:60073a}, Perkins \cite{MR88k:60085}, Jacka \cite{MR89j:60054}) then relaxed the assumption of finite second moment and presented constructions which for any centered probability measure $\mu$ give a stopping time $T$ such that $B_T\sim \mu$ and $(B_{t\land T}:t\ge 0)$ is a uniformly integrable martingale. We shall call stopping times for which the latter property is verified \emph{UI stopping times}.
These constructions work in the setting of continuous local martingales and some can be extended to specific discontinuous setups (cf.~Ob\l\'oj and Yor \cite{obloj_yor}).

When the target measure $\mu$ is not centered the process $(B_{t\land T}:t\ge 0)$ can not be a uniformly integrable martingale. Hence, a more general criterion for deciding when the stopping time $T$ is \emph{reasonably small} is needed
and such criterion is provided by notion of \emph{minimality} introduced by Monroe \cite{MR49:8096}, and considered more recently by Cox and Hobson \cite{cox_hobson2}. We say that a stopping time $T$ is \emph{minimal} if whenever $S\le T$ is a stopping time such that $B_S\sim B_T$ then $S=T$ a.s.. Imposing a minimality requirement on the solutions to the Skorokhod embedding problem is justified by a result of Monroe \cite{MR49:8096} which asserts that a stopping time $T$ which embeds a centered distribution
in a Brownian motion
is \emph{minimal} if and only if $T$ is a \emph{UI stopping time}. Recently Cox and Hobson \cite{cox_hobson2} and Cox \cite{cox_min} provided a description of minimality for general starting and target measures. Although it is not possible for the stopped process to be uniformly integrable in general, the conditions are closely related to uniform integrability in the Brownian setting. We can thus say that the notion of minimality for Brownian motion is well understood and is a feasible criterion.

Once we understand the equivalence between minimal and UI stopping times for Brownian motion (and via time-change arguments for all continuous local martingales) a natural question is then to ask: what is the situation for other martingales? More precisely, as we note below, uniform integrability always implies minimality, so the question is when, and `how much', is the former more restrictive?
In one sense, we would like to discover the `correct' interpretation of \emph{small} for embeddings, and determine when the different definitions agree and disagree. In general this question appears to be hard, and one of the aims of this work is to demonstrate that even in simple cases, it is not easy to classify the relevant sets. More precisely, in this work we will focus on the simple symmetric random walk, relative to its natural filtration. The latter assumption is important, since if we allowed for example enough external randomisation, we could reconstruct a Brownian motion from the random walk, and we would be returning to this setting. The restriction to the natural filtration will alter the problem sufficiently to provide interesting differences.

As well as the sets of minimal and UI embeddings, we also consider two natural constructions from the Brownian setting, and compare the distributions which may be embedded via these constructions. The constructions of interest here are the Azema-Yor \cite{MR82c:60073a} and Chacon-Walsh \cite{MR56:3934} embeddings; the latter of which can be considered as the set of stopping times which are the composition of first exit times from intervals.

The paper shall proceed as follows: in Section~\ref{sec:remarks} we make some initial remarks on the different classes of embeddings we shall consider; in Section~\ref{sec:potential} we examine the Azema-Yor and Chacon-Walsh constructions in the random walk setting, and Sections~\ref{sec:minimal} and \ref{sec:UIgeneral} look at UI embeddings in this context, where we find that the sets of admissible target measures can be remarkably complex.

\subsection*{Comments on notation}
In the sequel we deal mainly with processes in discrete time where time is indexed by $n=0,1,2,\dots$. When we refer to the continuous time setting time will be denoted by $t\in [0,\infty)$. Stopping times in discrete setting are denoted with Greek letters (typically $\tau$) and in continuous time with capital Latin letters (typically $T$). The set of probability measures on $\zr$ is denoted $\mathcal{M}$.\\
Throughout, $(X_n:n\ge 0)$ denotes a standard random walk, i.e.\ $X_n=\sum_{k=0}^n \xi_k$, where $(\xi_k)$ is a sequence of i.i.d.\ Bernoulli variables. The maximum is denoted by $\overline{X}_n=\max_{k\le n} X_k$
% and the \emph{local time in zero} is given as the number of visits to zero $l_n=\sum_{k=0}^n \mathbf{1}_{X_k=0}$
.\\
In the continuous time setup we will use $(B_t:t\ge 0)$ to denote a standard real-valued Brownian motion. Its maximum is $\overline{B}_t=\sup_{s\le t}B_s$.
% and its local time in zero as $L_t=\lim_{\epsilon\to 0}\frac{1}{2\epsilon}\int_0^t \mathbf{1}_{|B_s|\le \epsilon}ds$
A probability distribution is typically denoted by $\mu$. Its tail is given as $\overline{\mu}(x)=\mu([x,\infty))$. Dirac's delta at a point $x$ is denoted with $\delta_{x}$.

%\section{The Skorokhod embedding problem -- continuous time setting}
%
%

%
%The original motivation behind Skorokhod's work was to embed a random walk into a Brownian motion. We will use this technique also in our paper. Iterating any solution to the Skorokhod embedding problem we can obtain a sequence of stopping times $(T_n)$, $T_0=0$, such that $(B_{T_n}:n\ge 0)\sim(X_n:n\ge 0)$. We refer to Ob\l\'oj \cite[Sec.~11]{genealogia} for the details.
%
\section{The Skorokhod embedding for random walks: general remarks} \label{sec:remarks}
In this section we prove the existence of a minimal stopping time which solves the Skorokhod embedding problem for random walk, and make some simple observations which show that the discrete time setting is quite different from the continuous time setting. Of importance here is the fact that we are considering stopping times $\tau$ with respect to the natural filtration of the discrete process.

Under the assumption that we have additional information, we note that one can give a simple explicit randomised embedding which just requires an independent two-dimensional random variable. This can be done mimicking Hall's solution \cite{hall_techrep} (cf.~Ob\l\'oj \cite[Sec.~3.4]{genealogia}): for $\mu$ a centered probability distribution on $\zr$, $\sum_{k\ge 0}k\mu(\{k\})=m<\infty$, let $(U,V)$ be an independent variable with $\p(U=u,V=v)=\frac{(u-v)}{m}\mu(\{u\})\mu(\{v\})$, $u<0\le v$. Then $\tau=\inf\{n\ge 0: X_n\in \{U,V\}\}$ is a UI stopping time with $X_\tau \sim \mu$.

In fact, given suitable randomisation, we can even make a connection with the solutions for the Brownian case: given a random walk and sufficient independent randomisation, we are able to construct a Brownian motion by generating the intermediate paths, conditional on the start and end points, and further conditional on the end point being the first hitting time of a suitable integer. Now, given a stopping time $T$ for the Brownian motion, which embeds $\mu$ on $\Z$, we can construct a stopping time $\tau$ for the random walk (in an enlarged filtration) by considering the filtration $\td{\F}_n$ for the random walk generated by $(\F_{T_n}, \{T < T_{n+1}\})$ --- note that the martingale property ensures that $X_n$ remains a random walk in this filtration --- and defining $\tau$ by $\tau = n$ on $\{T_n \le T < T_{n+1}\}$. In particular, $X_{\tau} = B_T$ \as{}. It is clear that the stopping time $T$ is UI if and only if $\tau$ is since \eg{} $\sup_{t \le T} B_t$ and $\sup_{n \le \tau} X_n$ differ by at most 1.

Denote by $\mathcal{M}_0$ the set of all centered probability measures on $\mathbb{Z}$ and by
$\mathcal{M}^{UI}_0$ the set of probability measures $\mu$ on $\mathbb{Z}$ such that there exists a stopping time $\tau$ (in the natural filtration of $X$) such that $X_\tau \sim \mu$ and $(X_{n\land\tau}:n\ge 0)$ is a uniformly integrable martingale. Naturally, as the mean of a UI martingale is constant, we have $\mathcal{M}^{UI}_0\subseteq \mathcal{M}_0$. However, unlike in the setup of Brownian motion, the inclusion is strict:
\begin{proposition}
We have $\mathcal{M}_0^{UI}\subsetneq \mathcal{M}_0$.
\end{proposition}
\begin{proof}
To see this consider the target measure $\mu=\frac{1}{3}(\delta_{-1}+\delta_{0}+\delta_{1})$ which is a centered probability measure. Then if $\tau$ is a stopping time which embeds $\mu$, $X_\tau\sim \mu$ then $\tau$ can not be a UI stopping time, that is $(X_{n\land \tau}:n\ge 0)$ can not be a uniformly integrable martingale. This is simply because $\tau>1$, thus $(|X_{n\land \tau}|:n\ge 0)$ with positive probability goes above $1$ and does not stop until $|X_n|$ returns to $1$.\\
More precisely, if $\rho=\inf\{n>0:X_n=0\}$ then the process $(X_{n\land \rho}:n\ge 0)$ is not uniformly integrable and we have, $C>1$, $\e |X_{n\land \tau}|\mathbf{1}_{|X_{n\land \tau}|>C}\ge \e|X_{n\land \rho}|\mathbf{1}_{|X_{n\land \rho}|>C}$ so $(X_{n\land \tau}:n\ge 0)$ can not be uniformly integrable. Thus $\mu\in \mathcal{M}_0\setminus\mathcal{M}_0^{UI}$.
\end{proof}

It is a general fact, which holds for any real-valued martingale, that a UI embedding is minimal\footnote{The proofs in Monroe \cite[Thm~1]{MR49:8096} or Cox and Hobson \cite{cox_hobson2} even though written for Brownian motion generalise to an arbitrary martingale.}. The reverse is true in the Brownian motion setup with centered target laws, but not in general. It is thus natural to ask in the random walk setting: what measures can we embed in a minimal way?
The answer is given in the following theorem.
\begin{theorem}
\label{thm:existence_min}
For any probability measure $\mu$ on $\zr$ there exists a minimal stopping time $\tau$ with respect to the natural filtration of $(X_n)$ such that $X_\tau \sim \mu$.
\end{theorem}
\begin{proof}
Fix a probability distribution $\mu$ on $\zr$, which we write as $\mu=\sum_{k=0}^\infty a_k\delta_{k}$ with $a_k\ge 0$, $\sum a_k=1$. We first show that there exists a stopping time $\tau$ such that $X_\tau\sim \mu$. Let $Y_n=\Delta X_n=X_n-X_{n-1}$, $n\ge 1$, and recall that the random variable $U=\sum_{k\ge 1}2^{-k}\mathbf{1}_{Y_k=1}$ has a uniform distribution on $[0,1]$. Let $(\tilde{a}_i)$ be the decreasing reordering of the sequence $(a_k)$ and $(k_i)$ the corresponding reordering of atoms so that $\mu=\sum_{i=0}^\infty \tilde{a}_i\delta_{k_i}$, $\tilde{a}_i\ge \tilde{a}_{i+1}$ and $k_i\neq k_j$ for $i\neq j$. Write $\sigma$ for the reverse permutation: for a given $j$ we define $\sigma(j)$ via $\tilde{a}_{\sigma(j)}=a_j$.
Denote $\tilde{b}_i=\sum_{j=0}^i \tilde{a}_j$ and $N(U)$ the unique number $i$ such that $\tilde{b}_i\le U<\tilde{b}_{i+1}$. Note that $\p(N(U)=i)=\tilde{a}_i$. We claim that the value of $N(U)$ is a.s.\ determined after finite number of steps of the random walk $(X_n)$. Indeed, let $\rho=\inf\{n:Y_n=-1\}$ which is a.s.\ finite. Note that $U\in [0,1-2^{-\rho}]$ and there exists a random (finite) index $i_\rho$ such that $\tilde{b}_{i_\rho-1}\leq 1-2^{-\rho}<\tilde{b}_{i_\rho}$. Let $d=\min_{i\le i_\rho}|U-\tilde{b}_i|$. Then after $n_d:=(1-\log_2 d)$ steps of the random walk, as $U$ can differ at most by $d/2$ from $\sum_{k\ge 1}^{n_d}2^{-k}\mathbf{1}_{Y_k=1}$, $N(U)$ is uniquely determined. Finally, since $(X_n)$ is recurrent, the stopping time $\tau=\inf\{n>1-\log_2 d: X_n=k_{N(U)}\}$ is a.s.\ finite and it satisfies $X_\tau \sim \mu$ as $\p(X_\tau=j)=\p(k_{N(U)}=j)=\p(N(U)=\sigma(j))=\tilde{a}_{\sigma(j)}=a_j$.

It remains to see that we can choose $\tau$ minimal. This follows from standard reasoning (cf.\ Monroe \cite{MR49:8096}) as we know now that the set of stopping times $\{\tau: X_\tau \sim \mu\}$ is nonempty and it is naturally partially ordered (by $\preceq$, where $S \preceq T$ if and only if $S \le T$ \as{}.; see also Cox and Hobson \cite{cox_hobson2}).
\end{proof}
We can rewrite the theorem in short as $\mathcal{M}_0=\mathcal{M}^{MIN}_0$, where $\mathcal{M}_0^{MIN}$ denotes the set of centered probability measures on $\zr$ which can be embedded in a random walk by means of a minimal stopping time.

\section{Embeddings via potential theory} \label{sec:potential}
One-dimensional potential theory, as used by Chacon and Walsh \cite{MR56:3934}, proved a very useful tool for developing solutions to the Skorokhod embedding problem (cf.~Ob\l\'oj \cite{genealogia}). We apply it here in the framework of a random walk.
%\smallskip\\ \emph{In this section measure $\mu$ on $\zr$ is supposed integrable: $\sum_{n\in \zr}|n|\mu(\{n\})<\infty$.}\smallskip
In this section we suppose the measure $\mu$ on $\zr$ is integrable: $\sum_{n\in \zr}|n|\mu(\{n\})<\infty$.

Define the potential of $\mu$ on $\zr$ by
\begin{equation}\label{eq:potential}
    u_\mu(x)=-\int|x-y|d\mu(y)=-\sum_{n\in\zr}|x-n|\mu(\{n\}), \quad x\in \re.
\end{equation}
This is a continuous, piece-wise linear function breaking at atoms of $\mu$. We have $u_\mu(x)\leq -|x-\sum n\mu(\{n\})|$ with equality as $|x|\to \infty$. The potential function determines uniquely the measure and vice-versa. Furthermore, the pointwise convergence of potentials corresponds to the weak convergence of measures.
The crucial property for us lies in the fact that change in the potential of the distribution of a random walk is easy to characterise when stopped at first exit times. More precisely, let $\tau$ be a stopping time with $\e |X_\tau|<\infty$ and $\rho^\tau_{a,b}=\inf\{n\ge \tau: X_n\notin (a,b)\}$ for $a,b\in \zr$. Denote $u_1$ and $u_2$ the potentials of the distributions of $X_\tau$ and $X_{\rho^\tau_{a,b}}$ respectively. Then $u_2\le u_1$, $u_1(x)=u_2(x)$ for $x\notin (a,b)$ and $u_2$ is linear on $[a,b]$ . In other words, $u_2=\min\{u_1,l\}$ where $l$ is the line that goes through $(a,u_1(a))$ and $(b,u_1(b))$ (cf.~Chacon \cite{MR58:18746}, Cox \cite{cox_min}, Ob\l\'oj \cite[Sec.~2.2]{genealogia} for the details). We deduce the following fact.
\begin{lemma}
If there exists a sequence of linear functions $f_k$ with $|f_k'|<1$ such that $u_\mu=\lim u_k$, where $u_0(x)=-|x|$, $u_k=\min\{u_{k-1},f_k\}$ and $u_k$ is differentiable on $\re\setminus \zr$ then there exists a UI stopping time $\tau$ such that $X_\tau\sim \mu$.
\end{lemma}
%\noindent The condition that $u_k$ is differentiable on $\re\setminus \zr$ can be rephrased saying that $u_k$, which is point-wise linear, breaks only in points of $\zr$.
\begin{proof}
The conditions in the lemma imply $u_\mu\le u_0$ and thus $\mu$ is centered.
The stopping time $\tau$ is simply a superposition of first exit times. More precisely, consider a subsequence of $(u_k)$, which we still denote $(u_k)$, such that for every $k$ there exists $x_k$ such that $u_k(x_k)<u_{k-1}(x_k)$.
Define $a_k=\inf\{x:f_k<u_{k-1}\}$ and $b_k=\sup\{x: f_k< u_{k-1}\}$ and $\tau_k=\inf\{n\geq \tau_{k-1}:X_n\notin [a_k,b_k]\}$ with $\tau_0=0$. Note that with our assumptions, $a_k,b_k\in \zr$. Then $u_k$ is the potential of the law of $X_{\tau_k}$ and $\tau_k\nearrow \tau$ as $k\to \infty$. From the convergence of the potentials we deduce that $\tau$ is finite a.s.\ and $X_\tau \sim \mu$. The uniform integrability follows from standard arguments (cf.~Chacon \cite[Lemma 5.1]{MR58:18746}).
\end{proof}
We will call stopping times obtained in the above manner \emph{Chacon-Walsh stopping times} and the class of probability measures which can be embedded using these stopping times is denoted $\mathcal{M}_0^{CHW}$. We have $\mathcal{M}_0^{CHW}\subset \mathcal{M}_0^{UI}$ and the inclusion is strict. An example of an element of $\mathcal{M}_0^{UI}\setminus \mathcal{M}_0^{CHW}$ is given by $\mu=\frac{5}{16}\delta_{0}+\frac{11}{32}\delta_{-2}+\frac{11}{32}\delta_{2}$. That measure $\mu\in \mathcal{M}_0^{UI}$ will follow from Theorem \ref{thm:202char}. It is a tedious verification of all possibilities that $\mu\notin \mathcal{M}_0^{CHW}$, and is probably best seen graphically. It follows from the fact that  $u_\mu(0)=-\frac{11}{8}$, while when composing first exit times we can not have the value of the potential in $0$ in $(-\frac{3}{2},-\frac{4}{3})$. The value $-\frac{4}{3}$ is obtained via $\rho_{0,2}^{\rho^0_{-2,1}}$ and $-\frac{3}{2}$ via $\rho_{-1,1}^{\rho_{-2,0}^\tau}$, where $\tau=\rho_{0,2}^{\rho^0_{-1,1}}$.\\

Related to the Chacon-Walsh construction in the Brownian setting is the solution of Az\'ema and Yor \cite{MR82c:60073a}. For a centered probability measure $\mu$ on $\re$ define the \emph{Hardy-Littlewood} or \emph{barycenter} function via
\begin{equation}
\label{eq:barycenter}
    \Psi_\mu(x)=\frac{1}{\overline{\mu}(x)}\int_{[x,\infty)}y \, d\mu(y).
\end{equation}
Then the stopping time $T^\mu_{AY}=\inf\{t: \overline{B}_t\ge \Psi_\mu(B_t)\}$ embeds $\mu$, (\ie{} $B_{T^\mu_{AY}}\sim \mu$) and $(B_{t\land T^\mu_{AY}}:t\ge 0)$ is a uniformly integrable martingale.

With this in mind, we can consider a special case of the Chacon-Walsh construction in which the lines $f_n$ are tangential to $u_\mu$ and take them in a given order: from left to right. Then the sequences $(a_k)$ and $(b_k)$ are increasing and therefore $\tau$ is the first time we go below a certain level which is a function of the present maximum $\overline{X}$ (which basically tells us which of $b_k$ we have hit so far).
This corresponds to the solution of Az\'ema and Yor as observed by Meilijson \cite{MR86d:60052}.\footnote{The \emph{barycenter} function $\Psi_\mu(x)$ displayed in \eqref{eq:barycenter} can be seen as the intersection of the tangent to $u_\mu$ in point $x$ with the line $-|x|$ (cf.~Ob\l\'oj \cite[Sec.~5]{genealogia}).}
We have thus the following result\footnote{Similar remarks for discrete martingales were made in Fujita \cite{fujita} and Ob\l\'oj \cite[Sec.~4]{genealogia}.}.
\begin{proposition}
Let $\mu$ be a centered probability measure on $\zr$. The Az\'ema-Yor stopping time $T^\mu_{AY}=\inf\{n: \overline{X}_n\ge \Psi_\mu(X_n)\}$ embeds $\mu$ if and only if $\Psi_\mu$, displayed in \eqref{eq:barycenter}, satisfies $\Psi_\mu(x)\in \nr$. Then, $(X_{n\land T^\mu_{AY}}:n\ge 0)$ is a uniformly integrable martingale.
\end{proposition}
\begin{proof}
  Sufficiency of the condition was argued above. To see that it is also necessary recall (cf.\ Revuz and Yor \cite[p.~271]{MR2000h:60050}) the one to one correspondence, given by $\mu\to\Psi_\mu$, between centered probability measures $\mu$ on $\re$ and positive, left-continuous, non-decreasing functions $\Psi$ such that there exist $-\infty\le a<0<b\le \infty$, $\Psi(x)=0$ on $(-\infty,a]$, $\Psi(x)>x$ on $(a,b)$ and $\Psi(x)=x$ on $[b,\infty)$. Note that $\Psi_\mu$ is constant outside the support of $\mu$, so in particular when $\mu(\zr)=1$ then $\Psi_\mu$ is constant on every interval $(k,k+1]$. Then let $\mu$ be a probability measure on $\zr$ such that there exists $k\in \zr$ with $\Psi_\mu(k)\notin \nr$. Possibly choosing a different $k$ we can suppose that $\Psi_\mu(k)<\Psi_\mu(k+1)$ or equivalently that $\mu(\{k\})>0$. Let $l\in \nr$ be such that $l<\Psi_\mu(k)<l+1$. Then we either have $\Psi_\mu(k+1)\le l+1$ or $l+1< \Psi_\mu(k+1)$. In the first case the process will never stop in $k$, $\p(X_{T^\mu_{AY}}=k)=0$, which shows that $X_{T^\mu_{AY}}\nsim\mu$. In the second case, changing the value of $\Psi_\mu$ to any other value between $(\Psi_\mu(k-1)\lor l,l+1)$ will not affect the stopping time. We thus obtain a continuity of functions $\Psi$, each corresponding to a different measure on $\zr$, which all yield the same stopping time and thus the same law of the stopped process.
\end{proof}
We denote the class of measures which can be embedded using Az\'ema-Yor's stopping times with $\mathcal{M}_0^{AY}$. Naturally we have $\mathcal{M}_0^{AY}\subset \mathcal{M}_0^{CHW}$. Moreover, unlike in the continuous-time setup of Brownian motion, the inclusion is strict. To see this we recall an example given in Ob\l\'oj \cite[Sec.~4]{genealogia}: consider $\mu=\frac{2}{9}\delta_{-3}+\frac{4}{9}\delta_{0}+\frac{1}{3}\delta_{2}$. Then $\Psi_\mu(0)=\frac{6}{7}\notin \nr$. However the Chacon-Walsh stopping time $\inf\{n>\rho^0_{-1,2}: X_n\notin [-3,0]\}$, where $\rho^0_{-1,2}=\inf\{n\ge 0:X_n\notin [-1,2]\}$, embeds $\mu$.

Gathering the results described so far we conclude that
$$\mathcal{M}_0^{AY}\subsetneq \mathcal{M}_0^{CHW}\subsetneq \mathcal{M}_0^{UI}\subsetneq \mathcal{M}_0^{MIN}=\mathcal{M}_0$$
which is in sharp comparison with the continuous--time setup of Brownian motion\footnote{And therefore via time-change arguments, for any continuous local martingale, with \as{} infinite quadratic variation.} where all the sets are equal.

\section{UI embeddings: first steps}
\label{sec:minimal}
In this and the subsequent section we classify the possible elements of $\mathcal{M}_0^{UI}$.
Initially we consider measures with support on $\{-N, \ldots ,-1,0,1, \ldots ,N\}$. The restriction to $[-N,N]$ forces the candidate stopping times $\tau$ to satisfy $\tau \le \inf\{n \ge 0 : X_n \in \{-N,N\}\}$. As we shall see, requiring $\tau$ to be a stopping time in the natural filtration forces a complex, fractal structure on the set of possible hitting measures. For example, as a trivial initial statement, when $N=2$ we note that we cannot stop at zero with a probability in $(\frac{1}{2},1)$ --- either we stop at time $0$, with probability $1$, or else the first time we could stop will be time $2$, however with probability $\frac{1}{2}$ we will hit $\{-2,2\}$ before returning to 0.

We begin by concentrating on the case where the stopped distribution actually has support on the set $\{-2,0,2\}$. The analysis will depend on counting the number of possible paths after $2n$ steps. After $2n$ steps, there will be $2^{2n}$ possible paths, each occurring with equal probability, however only $2^n$ of these paths will not have hit $\{-2,2\}$, and all of these paths will be at $0$ at time $2n$. Since the stopping time $\tau$ is adapted to the natural filtration of $X$, if a path is stopped at $0$ at time $2n$, all paths which look identical up to time $2n$ must also stop at $2n$. Consequently, given a stopping time $\tau$, we can encode its properties in terms of the number of paths it will stop at time $2n$; we do this using the sequence $(a_0,a_1,a_2,\ldots)$, so that $a_n$ is the number of different (up to time $2n$) paths which are stopped by $\tau$ at time $2n$. We can also reverse the process, so that given a suitable sequence $(a_0, a_1, a_2, \ldots)$ we define a stopping time $\tau$ which stops at $0$ (according to some algorithm) $a_n$ different paths at time $2n$. Of course, not all sequences will necessarily allow such a stopping time to be defined, and the exact criteria are given in the following theorem.

\begin{theorem}
\label{thm:202char}
Let $\mu\in\mathcal{M}_0$ with support on $\{-2,0,2\}$, $\mu(\{0\}) = p=1-\mu(\{-2,2\})$. Then, $\mu \in \mathcal{M}_0^{UI}$  if and only if $p$ can be written as a base-4 fraction of the form $a_0.a_1a_2a_3\ldots$ with $a_n \in \{0,1,2,3\}$, where
\begin{equation} \label{eqn:anle}
a_n \le 2^n - \sum_{i=1}^n 2^i a_{n-i},
\end{equation}
or equivalently,
\begin{equation} \label{eqn:anle2}
\sum_{i \ge 0} 2^{-i} a_{i} \le 1.
\end{equation}
Furthermore, the set $\mathcal{S}$ of admissible values of $p=\mu(\{0\})$ is the unique fixed point of the mapping $s$ operating on the closed subsets of $[0,1]$ given by
\begin{equation}\label{eq:contr1}
    A\stackrel{s}{\to}[0,\frac{1}{8}]\cup
    \Big(\frac{1}{4}A+\frac{1}{8}\Big)\cup \Big(\frac{1}{4}A+\frac{1}{4}\Big)\cup\{1\}\ .
\end{equation}
\end{theorem}

\begin{proof}
Suppose that we have a probability $p = a_0.a_1 a_2\ldots$ satisfying \eqref{eqn:anle}; as remarked above, we can convert the sequence into a stopping time, however we must ensure that at each time $2n$, there exist sufficiently many different paths arriving to be able to stop $a_n$ paths. Suppose at time $2n$ there are $k_n$ paths, then we require $a_n \le k_n$. Assuming this is true, there will then be $2(k_n-a_n)$ different paths at $0$ at time $2(n+1)$, so by a similar reasoning, we must therefore have $a_{n+1} \le k_{n+1} = 2(k_n - a_n)$. Noting that $k_0 = 1$, we can iterate this procedure to deduce \eqref{eqn:anle}.

Conversely, given a stopping time $\tau$, we can derive a sequence $(a_0,a_1,\ldots)$ corresponding to the number of paths stopped at each stage. By the above argument, these $a_n$ satisfy \eqref{eqn:anle}; what is not necessarily true is that each $a_n \in \{0,1,2,3\}$. However the probability of stopping at $0$ is still given by $\sum_{i\ge 0} 4^{-n} a_n$, and we can form a new sequence $(\tilde{a}_0,\tilde{a}_1,\tilde{a}_2,\ldots)$ such that $\tilde{a}_n \in \{0,1,2,3\}$ and $\sum_{i\ge 0} 4^{-n} \tilde{a}_n = \sum_{i\ge 0} 4^{-n} a_n$. Where necessary we will work with a sequence which terminates in a string of zeros rather than a string of threes. However for such a sequence, it is then clear that
\[
\sum_{i=0}^\infty 2^{-i} \tilde{a}_i \le \sum_{i=0}^\infty 2^{-i} a_i
\]
(replacing a $4$ in the $j^{th}$ position with a $1$ in the $(j-1)^{th}$ position always reduces the value, and the total value of the sum is bounded above by $1$, and below by $0$), so that the result holds in general.\smallskip\\
It remains to prove the last assertion of the theorem. Define
set functions, mapping the set of closed subsets of $[0,1]$ to itself via, $A\subset [0,1]$,
\renewcommand{\arraystretch}{1.5}
\begin{equation}
\label{eqn:system1}
\begin{array}{rclcrclcrcl}
f_1(A) & = & \frac{1}{4} + \frac{1}{4}A & \ & f_2(A) & = & \frac{1}{8} + \frac{1}{4}A & \ &
f_3(A) & = & \left[0,\frac{1}{8}\right] \cup\{1\}
\end{array}
\end{equation}
For convenience, when dealing with singletons $\{p\}$ we write simply $f_1(p)=1/4+1/4p$ etc. Note that $s(A)=f_1(A)\cup f_2(A)\cup f_3(A)$. It is now clear from the definition of $s$ that it is a contraction mapping under the Hausdorff metric\footnote{If $X$ is a metric space, the Hausdorff metric is defined on set of compact subsets $A,B$ of $X$ by
\[
d_H(A,B) = \inf \{r > 0 : d(A,y) \le r\ \forall y \in B \mbox{ and } d(x,B) \le r\ \forall x \in A\}.
\]
}, and hence, by the Contraction Mapping Theorem, has a unique fixed point in the set of compact subsets of $\re$. It is simple to check that $\mathcal{S}$ is a closed subset of $[0,1]$ (by considering for example the base-4 expansions), thus our goal is to show that $s(\mathcal{S})=\mathcal{S}$.

We first show that $s(\mathcal{S})\subset\mathcal{S}$. To see this we simply check that if $p\in\mathcal{S}$ then $f_i(p)\in \mathcal{S}$ for $i=1,2$ and that $[0,\frac{1}{8}]\cup\{1\}\subset \mathcal{S}$. Consider for example $f_1$. The case $p=1$ is trivial. Let $p\in \mathcal{S}$, $p<1$, and write it in base-4 expansion as $0.a_1a_2\dots$. Then $f_1(p)=0.1a_1a_2\dots$ and \eqref{eqn:anle2} holds: so by the first part of the theorem $f_1(p)\in \mathcal{S}$. We proceed likewise for $f_2$. Finally, to prove $[0,\frac{1}{8}]\subset\mathcal{S}$, take any $0<p<\frac{1}{8}$ and write its base-4 expansion $p=0.0a_2a_3\dots$ where $a_2\in \{0,1\}$. Then $\sum_{i=0}^\infty a_i2^{-i}\leq \frac{1}{4}+3\sum_{i=3}^\infty 2^{-i}=1$ which shows that $p\in \mathcal{S}$.\\
It remains to see the converse, namely that $\mathcal{S}\subset s(\mathcal{S})$. Let $p\in \mathcal{S}$ and write its base-4 expansion $p=a_0.a_1 a_2 a_3\dots$. We will analyse various cases and use implicitly the criterion \eqref{eqn:anle2}.
The case $p=1$ is trivial we can therefore suppose $a_0=0$. If $a_1=2$ then $p=1/2$ and we have $p=f_1(1)$. If $a_1=1$ then $p=f_1(q)$ with $q=0.a_2 a_3 a_4\dots$. To see that $q\in \mathcal{S}$ note that since $p\in \mathcal{S}$ we have $1/2+\sum_{i=2}^\infty 2^{-i}a_i \le 1$ and thus $\sum_{i=1}^\infty 2^{-i}a_{i+1}\le 1$.\\
Suppose now that $a_1=0$. If $a_2=3$ then $p=f_2(q)$ with $q=0.1a_3 a_4\dots$ and again since $p\in \mathcal{S}$ we have $\sum_{i=3}^\infty 2^{-i}a_i\le 1/4$ which implies that $q\in \mathcal{S}$. If $a_2=2$ then $p=f_2(q)$ with $q=0.0a_3 a_4\dots$ and we check again that $q\in \mathcal{S}$. Finally if $a_2\le 1$ then $p<1/8$ and is thus in the image of $f_4$.
We obtain finally that $f(\mathcal{S})=\mathcal{S}$ and thus $\mathcal{S}$ is the fixed point of the contraction mapping $s$ which ends the proof of the theorem.
\end{proof}

We want to comment the rather surprising nature of the set $\mathcal{S}$. It is in fact a self-similar structure, or fractal. In particular, following the characterisation of \cite{MR799111} (see also Falconer \cite[Chap.~9]{MR2118797}), we can say that $\mathcal{S}$ is an {\it iterated function system with a condensation set} generated by the system \eqref{eqn:system1}. From the representation \eqref{eq:contr1} it is easy to deduce that the one-dimensional Lebesgue measure of $\mathcal{S}$ is equal to $\frac{1}{4}$.\\
An alternative representation of the set can also be given in which the set is the fixed point of a standard \emph{iterated function system}; that is, we can drop the condensation set, in exchange for a larger set of functions. We replace the function $f_3$ by function(s) $g$ which map $\mathcal{S}$ into $\mathcal{S}$ and $[0,1/8]$ onto $[0,1/8]$. To this end define $g_i(x)=\frac{1}{4}x+\frac{k}{64}$. Note that $g_{8}=f_2$ and $g_{16}=f_1$. We claim that the set $\mathcal{S}$ is the unique fixed point of the mapping
\begin{equation}\label{eqn:contr2}
    A\stackrel{\tilde{s}}{\to} g_0(A)\cup g_2(A)\cup g_4(A)\cup g_6(A)\cup g_8(A)\cup g_{16}(A)\cup \{1\}\ .
\end{equation}
It is immediate that $[0,1/8]\subset \tilde{s}([0,1/8])$. It remains to see that if $p\in \mathcal{S}$ then $g_i(p)\in \mathcal{S}$ for $i=0,2,4,6$ which is easily checked with \eqref{eqn:anle2}.\\
To deduce some more information about the structure of $\mathcal{S}$, observe that $g_8([0,1/8])=[1/8,5/32]$. Iterating this we see that $[0,x_*]\subset \mathcal{S}$ where $x_*$ satisfies $x_*=1/8+1/4 x_*$. We have thus $x_*=1/6$ which has $0.022222\dots$ base-4 expansion and corresponds to stopping $2$ trajectories every second step of the random walk starting with the $4^{th}$ step.\\
%Another natural question to ask concerns the multifractal spectrum of the set, however it is not hard to check, using for example the techniques of Lau and Ngai \cite{MR1667146}, that the multifractal spectrum of the resulting set is not interesting, and is concentrated at 1.\bigskip\\
%Consider now the more general case, where we allow the process to be stopped at $\{-1,1\}$ in addition. A similar analysis to before allows us to classify the measures $\mu \in \mathcal{M}_0^{UI}$ with support on $\{-2,-1,0,1,2\}$.
Another natural question to ask concerns the dimension of the set. It is clear that the presence of the interval $[0,1/8]$ forces the dimension of the whole set to be 1, however is this also true locally? It turns out that the local dimension of any point in the set is either 0 or 1. This can be seen relatively easily: consider a point $x \in \mathcal{S}$; either the base-4 expansion of this point is terminating (that is, can be written with a finite number of non-zero $a_n$) or it is not. In the latter case, given $r>0$, we can find $n$ such that $4^{-n} \le r < 4^{-n+1}$. Since the sequence we choose is not terminating, the value $k_{n+2}$ defined in the previous theorem is at least $1$; further, by defining a new set of points which agree with $x$ up to $a_{n+1}$, and have $a_{n+2} = a_{n+3}= 0$ we may take any other terminating sequence beyond this point. This interval of points therefore has Lebesgue measure at least $4^{-n-4}$, and is also contained in the ball of radius $4^{-n-1}$ about $x$. More specifically, (writing $B(x,r)$ for the ball with centre $x$ and radius $r$) we have $|B(x,r) \cap \mathcal{S}| \ge r 4^{-4}$ and
\[
\liminf_{r \to 0} \frac{\log(|B(x,r) \cap \mathcal{S}|)}{\log{r}} \ge 1.
\]
Since our set is a subset of $\R$, it is clear that the local dimension cannot exceed one at a non-terminating point. In the second case, consider a terminating point --- here there are two possibilities: either the $k_n$s are zero for sufficiently large $n$, in which case the point is isolated (there is clearly a small interval above the point which is empty, and it can similarly be checked that there is a small interval below the point), or the $k_n$'s increase after the final non-zero $a_n$, but in this case it is clear that there is a small interval of points above $x$, and as claimed, the point is either isolated, or has a local dimension of 1.\footnote{According to some definitions, the set we have described would not be a fractal, in that it has no non-integer dimensions even at the local level; however we follow the more general classification described in the introduction to Falconer \cite{MR2118797}, and note that the set clearly has a complex local structure, and exhibits many of the features typical of the more restrictive definition.}

\begin{theorem}
\label{thm:202fullchar}
Suppose that $\mu \in \mathcal{M}_0$ with support on $\{-2,-1,0,1,2\}$. Then
$\mu \in \mathcal{M}_0^{UI}$ if and only if
\begin{equation}
\label{eqn:3atomsvalues}
\mu(\{0\})=\sum_{i \ge 0} a_i 2^{-2i},\ \mu(\{-1\})=\sum_{i \ge 0} b_i 2^{-2i+1},\ \mu(\{1\})=\sum_{i \ge 0} c_i 2^{-2i+1}
\end{equation}
where $a_i,b_i,c_i\in\{0,1,2,3\}$ and the sequences satisfy: $b_0=c_0=0$,
\begin{eqnarray}
\sum_{i=0}^\infty 2^{-i}a_i+\sum_{i=0}^\infty 2^{-i}(b_i+c_i)&\le & 1\label{eqn:3atoms1cond}\\
2^{n-1} - \sum_{i =0}^{n-1} 2^{n-i-1} a_i - \sum_{i=0}^{n-1}2^{n-i-1}(b_i + c_i)&\ge & b_n\lor c_n,\ n\ge 0\ .\label{eqn:3atoms3cond}
\end{eqnarray}
Furthermore, the set $\mathcal{S}^{(3)}$ of possible values of $p=(\mu(\{-1\},\mu(\{0\}),\mu(\{1\}))$ is the unique fixed point of the mapping $f$ operating on the closed subsets of $[0,1]^3$, given by $A\mapsto \bigcup_{q\in \mathcal{Q}}(\frac{1}{4}A+q)\cup g(A)$, where $\mathcal{Q}$ is a finite set to be described in the proof and $g(A)=\{(0,1,0),(\frac{1}{2},0,\frac{1}{2}),(\frac{1}{2},\frac{1}{4},0),(0,\frac{1}{4},\frac{1}{2})\}$.
\end{theorem}
\begin{proof}
We have a picture similar to the one described before Theorem \ref{thm:202char}. As before, our approach will be to count the number of `different' paths, however we now need to consider stopping at all the points $-1,0,1$, and the corresponding constraints on the system. As before, $a_n$ will be identified with the number of paths which are stopped at $0$ after $2n$ steps, and we also now introduce the sequences $(b_{n})_{n \ge 1}$ and $(c_{n})_{n \ge 1}$ which will correspond to the stopping behaviour, after $(2n-1)$ steps, at the points $-1$ and $1$ respectively.\\
%We can think about a superposition of rhombuses where the central points correspond to the visits in zero and the middle points to visits in $-1$ and $1$. Every time a path visits zero it splits into two paths - one travelling down to $-1$ and one travelling up to $1$. After which they meet again in $0$, each of them splits in two etc.
%In each vertexes of $n^{th}$ rhombus we place values $a_{n-1}$ in the upper vertex, $b_n$ and $c_n$ in the left and right vertex respectively and $a_n$ in the lower vertex. See Figure for a graphical representation.\\ {\bf to do: figure here!!!}\\ These values correspond to number of paths we stop. Thus, for example, $b_n$ is the number of paths we stop in $-1$ after $2n-1$ steps.
%
% Is a figure really necessary here? Surely it is just the same as the previous case...
%
As before, we can also identify $p=(p_{-1},p_0,p_1)\in\mathcal{S}^{(3)}$ with the sequence $(b^{(p)}_n,a^{(p)}_n,c^{(p)}_n)_{n\ge 0}$, which is the base-4 expansion of $(\frac{p_{-1}}{2},p_0,\frac{p_1}{2})$. This we can transform into a stopping time provided that there are always enough paths to stop the prescribed number at each step. Denote $k_n^{(p)}$ the number of paths still arriving at $0$ after $2n$ steps, where in the first $(2n-1)$ steps we were successfully realizing the stopping rule prescribed by $p$. We drop the superscript $(p)$ when $p$ is fixed. Then we have to require that $a_n\le k_n$ and $b_n\le k_n-a_n$, $c_n\le k_n-a_n$. Using induction we can prove that
\begin{equation}
k_n  =  2^n - \sum_{i =0}^{n-1} 2^{n-i} a_i - \sum_{i=1}^{n}2^{n-i}(b_i + c_i)\label{eqn:3atoms2cond} .
\end{equation}
Then the condition $a_n\le k_n$, for all $n\ge 0$, can be rewritten under equivalent form \eqref{eqn:3atoms1cond}.
Note that it also contains the necessary condition on $(b_n+c_n)$, namely that $b_n+c_n\le 2(k_n-a_n)$. However, \eqref{eqn:3atoms1cond} does not encode the restriction $b_n\lor c_n\le k_n-a_n$, which is \eqref{eqn:3atoms3cond}.

Conversely, given a UI stopping time $\tau$ with $X_\tau\in \{-2,-1,0,1,2\}$ we can derive the sequence $(b_n,a_n,c_n)$ of paths stopped respectively in $(-1,0,1)$ after $(2n+1,2n,2n+1)$ steps. By the arguments above $(b_n,a_n,c_n)$ satisfy \eqref{eqn:3atomsvalues}, \eqref{eqn:3atoms1cond} and \eqref{eqn:3atoms3cond} but it is not necessarily true that $a_n,b_n,c_n\in \{0,1,2,3\}$. Suppose then that the sequence $(b_n,a_n,c_n)$ is terminating (\ie{} there exists $n_0$ such that $a_n=b_n=c_n=0$ for $n \ge n_0$), and for some $j\ge 3$ we have $b_j\ge 4$. Define a new sequence $(\tilde{b}_n,a_n,\tilde{c}_n)$ via $\tilde{b}_n=b_n$, $\tilde{c}_n=c_n$ for $n$ different from $j$ and $(j-1)$, and with $\tilde{b}_{j-1}=b_{j-1}+1$, $\tilde{b}_j=b_j-4$ and likewise $\tilde{c}_{j-1}=c_{j-1}+\mathbf{1}_{c_j\ge 4}$, $\tilde{c}_j=c_j-4\mathbf{1}_{c_j\ge 4}$. It is easy to see that the new sequence satisfies \eqref{eqn:3atoms1cond} and \eqref{eqn:3atoms3cond}. It thus encodes a stopping time $\tilde{\tau}$ and by \eqref{eqn:3atomsvalues} $X_\tau\sim X_{\tilde{\tau}}$. Iterating this argument we can assume that the sequence $(\tilde{b}_n,a_n,\tilde{c}_n)$ satisfies $\tilde{b}_n,\tilde{c}_n\in \{0,1,2,3\}$. Now taking $(\tilde{a}_n)$ as the base-4 expansion of $\sum_{n\ge 0}a_n 4^{-n}$ we have $\tilde{a}_n,\tilde{b}_n,\tilde{c}_n\in \{0,1,2,3\}$ and we verify immediately that $(\tilde{a}_n,\tilde{b}_n,\tilde{c}_n)$ satisfies \eqref{eqn:3atoms1cond}, \eqref{eqn:3atoms3cond} and encodes the same measure as $(b_n,a_n,c_n)$.

It remains now to show that the same can be said for a general sequence $(b_n,a_n,c_n)$. Let $p=(p_{-1},p_0,p_{1})$ be the associated point in $\mathcal{S}^{(3)}$ and $(\td{b}_n,\td{a}_n,\td{c}_n)$ its base-4 expansion. As in the proof of Theorem~\ref{thm:202char}, the latter satisfies \eqref{eqn:3atoms1cond} so all we need to show is that it also satisfies \eqref{eqn:3atoms3cond}.\\
First note that as \eqref{eqn:3atoms1cond}--\eqref{eqn:3atoms3cond} hold for $(b_n, a_n, c_n)$, they will also hold for the truncated sequences $(b_n^j,a_n^j,c_n^j)$, where the $j$ denotes $a^j_{j+1} = b^j_{j+1}=c^j_{j+1}= a^j_{j+2}= \ldots = 0$, and therefore, by the argument above, also for their base-4 expansions $(\td{b}_n^j,\td{a}_n^j,\td{c}_n^j)$. Observe that when two expansions exists we take the finite one. We will now argue that for any fixed $m$, for $j$ big enough, the sequences   $(\td{b}_n^j,\td{a}_n^j,\td{c}_n^j)$ and  $(\td{b}_n,\td{a}_n,\td{c}_n)$ coincide for $n\le m$, which will imply that the last sequence also satisfies \eqref{eqn:3atoms3cond}.\\
More precisely, we need to show that
\begin{equation}\label{eqn:proof4base}
    \forall m\; \exists j_m\; \forall j\ge j_m, \;(\td{b}_n^j,\td{a}_n^j,\td{c}_n^j)=(\td{b}_n,\td{a}_n,\td{c}_n)\textrm{ for } n<m.
\end{equation}
The argument is the same for all three sequences, so we present it for the sequence $(b_n)$. If it was terminating then obviously for $j$ larger than its length $(\td{b}_n)=(\td{b}^j_n)$. Suppose $(b_n)$ is not terminating. Note that $(b_n)$ is not terminating if and only if $(\td{b}_n)$ is not terminating. Let $p_{-1}^j=2\sum_{i=1}^j4^{-i}\td{b}^j_i$. Since we have also $p_{-1}^j=2\sum_{i=1}^j4^{-i}b_i$, we know that $p_{-1}^j\nearrow p_{-1}$ as $j\to \infty$. Fix $m> 1$ and let $q_m=\frac{p_{-1}}{2}-\sum_{i=1}^m4^{-i}\td{b}_i$.
% which is the distance of $\frac{p_{-1}}{2}$ from the nearest point on $4^{-m}$ grid to its left.
Then there exists $j_m$ such that for all $j\ge j_m$,  $p_{-1}-p_{-1}^j<2q_m$, which we can rewrite as $\sum_{i=1}^m4^{-i}\td{b}_i< \frac{p_{-1}^j}{2}\le \frac{p_{-1}}{2}$. The last inequality together with the obvious inequality $\frac{p_{-1}}{2}<\sum_{i=1}^m\td{b}_i+4^{-m}$, imply that base-4 expansions of $\frac{p_{-1}^j}{2}$ and of $\frac{p_{-1}}{2}$ coincide up to $m^{th}$ place, that is $\td{b}_i=\td{b}^j_i$ for all $i\le m$. The same argument applies to $(\td{a}_n)$ and $(\td{c}_n)$. This proves \eqref{eqn:proof4base} and consequently that the sequence $(\td{a}_n,\td{b}_n,\td{c}_n)$ satisfies \eqref{eqn:3atoms3cond}, which ends the proof of the first part of the theorem.\smallskip

We now move to the second part of the theorem. We could do an analysis as in Theorem \ref{thm:202char} however this would be very involved in the present setup. Instead, we generalise the technique used to arrive at \eqref{eqn:contr2}; as a consequence, we do not have a neat description of the functions, but rather an algorithm for obtaining them.\\
% We will rather present a general approach and it will be clear that in fact proving Theorem \ref{thm:202char} we also followed this approach.
The following observation proves to be crucial: if some $k_n^{(p)}$ is large enough then any sequence of $(a_i,b_i,c_i)_{i\ge n}$ is admissible. More precisely as $a_n,b_n,c_n\le 3$ we have $k_{n+1}^{(p)}= 2(k_n^{(p)}-a_n)-(b_n+c_n)\ge 2k_n^{(p)}-12$ and thus if at some point $k^{(p)}_n\ge 12$ then for all $m\ge n$ $k_m^{(p)}\ge 12$.\\
As the first consequence note that $k^{(0)}_4=16$ and thus any $p\in [0,1]^3$ such that $a_i^{(p)}=b_i^{(p)}=c_i^{(p)}=b_4^{(p)}=c_4^{(p)}=0$ for $i=0,1,2,3$ is in fact an element of $\mathcal{S}^{(3)}$.\\
Define $\mathcal{Q}$ as the set of all $q\in [0,1]^3$ such that $a_5^{(q)}=a_i^{(q)}=b_i^{(q)}=c_i^{(q)}=0$ for all $i> 5$ and $k_5^{(q)}\geq 16$.
$\mathcal{Q}$ is thus the set of probabilities which encode stopping strategies for the first $9$ steps of the random walk and which stop at most $16$ out of $32$ paths which come back to zero after $10$ steps. This is a finite set (its cardinality is trivially smaller then $4^{14}$ and is actually much smaller).
Denote $f_q(p)=p/4+q$. Note that for any $p\in \mathcal{S}^{(3)}$, $k_5^{(p/4)}=16+k_4^{(p)}\ge 16$ so that $f_q(p)=p/4+q\in\mathcal{S}^{(3)}$ for any $q\in \mathcal{Q}$.
This shows that $f(\mathcal{S}^{(3)})\subset \mathcal{S}^{(3)}$.\\
Conversely, take any $p\in\mathcal{S}^{(3)}$ with $p\notin \{(0,1,0),(\frac{1}{2},0,\frac{1}{2}),(\frac{1}{2},\frac{1}{4},0),(0,\frac{1}{4},\frac{1}{2})\}$ as these values (extremal points) are by definition in $f(\mathcal{S}^{(3)})$. If $b_1^{(p)}=1$ then $p=f_{(1/2,0,0)}(w)$ where $a^{(w)}_n=a_{n+1}^{(p)}$, $n\ge 0$, and $c^{(w)}_0=b^{(w)}_0=0$, $c^{(w)}_n=c_{n+1}^{(p)}$, $b^{(w)}_n=b_{n+1}^{(p)}$ for $n\ge 1$. Likewise, if $c_1^{(p)}=1$ then $p\in f_{(0,0,1/2)}(\mathcal{S}^{(3)})$. Finally, if $a_1^{(p)}=2$ then $p=f_{0}((0,1,0))$ and if $a_1^{(p)}=1$ then $p\in f_{(0,1/4,0)}(\mathcal{S}^{(3)})$.\\
We can therefore assume that $b_1^{{\scriptstyle (p)}}=c_1^{(p)}=a_0^{(p)}=a_1^{(p)}=0$ and present the general
argument. We will reason according to the value of $k_5^{(p)}$.
Suppose that $k_5^{(q)}\ge 16$, which means that the stopping strategy encoded by $p$  stops (in the first $9$ steps of the random walk) less than $16$ out of the $32$ paths which come back to zero after $10$ steps. Thus `this part' of $p$ is an element of $\mathcal{Q}$: put $q=(2\sum_{i=0}^5b^{(p)}_i/2^i,\sum_{i=0}^4a^{(p)}_i/2^i,2\sum_{i=0}^5c^{(p)}_i/2^i)$ then $q\in \mathcal{Q}$. Furthermore, $k_5^{(p-q)}=32$ and thus $k_4^{(4(p-q))}=16$ which as we know is enough to support any sequence of $(a_n,b_n,c_n)$ onwards. Thus $p\in f_q(\mathcal{S}^{(3)})$.
Finally, suppose that $k_5^{(p)}<16$, that is $p$ stops (in the first $9$ steps of the random walk) more than $16$ out of the $32$ paths which come back to zero after $10$ steps. Then there exists a $q\in \mathcal{Q}$ (possibly many of them) which encodes the way $p$ stops $16$ paths. That is, there exists $q\in \mathcal{Q}$ such that $k_5^{(p-q)}=k_5^{(p)}+16$ and thus $p\in f_q(\mathcal{S}^{(3)})$.
\end{proof}
The set $\mathcal{Q}$ arising in the proof would appear to be rather large. A careful analysis could probably bring down its size considerable yielding a significantly smaller iterated function set describing $\mathcal{S}^{(3)}$. We note that the possible values of $\mu(\{0\})$ are not changed. Put differently
\[
 \mathcal{S}^{(3)}\cap \left(\{0\}\times [0,1]\times\{0\}\right)=\mathcal{S}.
\]

\section{UI embeddings: general study}
\label{sec:UIgeneral}

We now turn to the analysis of arbitrary $\mu\in\mathcal{M}^{UI}_0$. Initially we consider the measures with finite support, and show that by taking suitable closures, we can classify the whole set.\\ Fix $N>1$.
Let $\mathcal{S}^{(2N+1)}\subset [0,1]^{2N+1}$ denote the set of probability measures $\mu\in \mathcal{M}_0^{UI}$ with support in $[-(N+1),N+1]$. More precisely $p\in \mathcal{S}^{(2N+1)}$, $p=(p_{-N},\dots,p_N)$ defines uniquely a centered probability measure $\mu_p$ with $\mu_p(\{i\})=p_i$, $|i|\le N$, $\mu_p(\{-(N+1),N,\dots,N,N+1\})=1$.\\
Let $(a^i_n)_{-N\le i\le N,n\ge 0}$ be an infinite matrix of integers. Its entries will correspond to number of stopped paths: $a^{2i+1}_n,a^{2i}_n$ will represent number of paths stopped respectively in $(2i+1)$ after $(2n-1)$ steps in $2i$ after $2n$ steps. With respect to the notation used in Theorem \ref{thm:202fullchar} we have $b_n=a^{-1}_n$ and $c_n=a^1_n$.
Define the matrix $(k^i_n)_{-N\le i\le N,n\ge 0}$ via
\begin{equation}
\label{eqn:ksystem}
\left\{
\begin{array}{lcl}
    k^i_0 &=&\mathbf{1}_{i=0},\\
    k^{2i+1}_{n+1}&=& k^{2i}_{n}-a^{2i}_{n}+k^{2(i+1)}_{n}-a^{2(i+1)}_{n},\ n\ge 0,\ 2i+1\in [-N,N],\\
    k^{2i}_{n+1}&=&
k^{2i+1}_{n+1}-a^{2i+1}_{n+1}+k^{2i-1}_{n+1}-a^{2i-1}_{n+1},\ n\ge 0,\ 2i\in [-N,N],\\
\end{array}%
\right.
\end{equation}
where we put $a^{N+1}_n=k^{N+1}_n$ and $a^{-(N+1)}_n=k^{-(N+1)}_n$. We think of $k^{2i}_n$ (resp. $k^{2i+1}_n$) as the number of paths arriving at $2i$ (resp. $2i+1$) after $2n$ (resp. $2n-1$) steps. We note that if all $a_n^i=0$ then $(k^i_n)$, $n\le N/2$, form the first $N$ rows of Pascal's triangle.
\begin{theorem}
Let $\mu\in\mathcal{M}_0$ with support in $\{-(N+1),\ldots,(N+1)\}$. Then $\mu\in \mathcal{M}_0^{UI}$ if and only if there exists a matrix of integers $(a^i_n)_{-N\le i\le N,n\ge 0}$ such that
\begin{equation}\label{eqn:genatoms}
\mu({i})=2^{(i\bmod{2})}\sum_{j=0}^\infty 4^{-j}a^{i}_j\quad\textrm{and}\quad a_n^{i}\le k^{i}_n,
\end{equation}
$i\in [-N,N]$, $n\ge 0$, where $(k^i_n)_{-N\le i\le N,n\ge 0}$ is defined via \eqref{eqn:ksystem}.\\
Furthermore, the set $\mathcal{S}^{(2N+1)}$ of such measures $\mu$ is the unique fixed point of the mapping $f$ operating on the closed subsets of $[0,1]^{(2N+1)}$, given by $A\mapsto \bigcup_{q\in \mathcal{W}}(\frac{1}{4}A+q)\cup g(A)$, where $\mathcal{W}$ is a compact set to be described in the proof and $g(A)=\{{\scriptstyle (0,\ldots,0,1,0,\dots, 0),(0,\ldots,\frac{1}{2},0,\frac{1}{2},\ldots, 0),(0,\ldots,\frac{1}{2},\frac{1}{4},0,\ldots,0),(0,\ldots,0,\frac{1}{4},\frac{1}{2},\ldots,0)}\}$.
\end{theorem}
\noindent \textbf{Remarks:}\smallskip\\
%\begin{remarks}
%\begin{itemize}
%\item
The most surprising aspect of this theorem is the second part which shows that for any $N$ the set $\mathcal{S}^{(2N+1)}$ has a complex self-similar structure.\smallskip\\
%\item
Note that we do not present any canonical manner to associate a unique matrix $(a^i_n)$ to a given $\mu$. This due to the fact that, in contrast with the results of Section~\ref{sec:minimal}, we can not assume that $a^i_n\in \{0,1,2,3\}$. To convince herself, we invite the reader to consider the measure $\mu=\frac{3}{4}\delta_0+\frac{1}{8}(\delta_{-4}+\delta_4)$ which has the associated (unique) matrix $(a^i_n)$ given by $a^i_n=0$ for $i\neq 0$ and $a^0_0=0$, $a^0_1=2$, $a^0_n=2^{n-1}$, $n\ge 2$, and which encodes the stopping time $\inf\{n>0:X_n\in\{-4,0,4\}\}$.
%\item
\smallskip\\
We observe that equations \eqref{eqn:ksystem} and \eqref{eqn:genatoms} are not in a closed form as before but rather have a recursive structure. Possibly a closed form may be derived but for practical verification and implementation the recursive form seems more suitable.
%\item
\smallskip\\
We can consider stopping times which stop maximally 3 paths in a given point at a given step. Then the reasoning presented in the proof of Theorem \ref{thm:202fullchar} applies: it suffices to ensure that at least $12$ paths arrive in a given point to secure feasibility of any subsequent stopping strategy in that point. We see thus that (suppose $N\ge 3$) any point $p$ with $p_i\le 4^{-|i|-1}\land 4^{-3}$ belongs to $\mathcal{S}^{(2N+1)}$. In particular, $\mathcal{S}^{(2N+1)}$ has positive $(2N+1)$-dimensional Lebesgue measure.
%\end{itemize}
%\end{remarks}
\begin{proof}
The theorem is a generalised version of our earlier detailed studies presented in Theorems \ref{thm:202char} and \ref{thm:202fullchar}. The first part of the theorem follows from our description of possible stopping times in the natural filtration of $(X_n)$. Integers $(a^i_n)$ and $(k^i_n)$ have the interpretation indicated above and the condition $a^i_n\le k^i_n$ ensures that there are enough paths arriving at $i$ after $2n$ ($2n-1$ for $i$ odd) steps to realise the prescribed stopping strategy. Note that in particular, as $a^i_n\ge 0$ and $k^i_0=0$ for $i\neq 0$ we have that $a^i_n=k^i_n=0$ for $n<i/2$.\\
There are two paths which come back to zero after 2 steps. Define $\mathcal{W}$ as the set of these points in $\mathcal{S}^{(2N+1)}$ which never stop descendants of at least one of these two paths:
$\mathcal{W}=\{p\in \mathcal{S}^{(2N+1)}: p+{\scriptstyle (0,\ldots,0,\frac{1}{4},0,\ldots, 0)}\in \mathcal{S}^{(2N+1)}\}$. The difference with the set $\mathcal{Q}$ defined in the proof of Theorem \ref{thm:202fullchar} is that there we considered only $p$ with base-4 expansions terminating after $5$ digits.
Observe that for any $p\in \mathcal{S}^{(2N+1)}$ and $q\in \mathcal{W}$, $f_q(p)=p/4+q\in \mathcal{S}^{(2N+1)}$ (this is simply because one path originating from zero after the second step suffices to ensure the stopping strategy prescribed by $p/4$). Conversely, for any $p\in \mathcal{S}^{(2N+1)}\setminus g(\mathcal{S}^{(2N+1)})$
we can find $q=q(p)\in \mathcal{W}$ such that $p\in f_q(\mathcal{S}^{(2N+1)})$ that is $4(p-q)\in \mathcal{S}^{(2N+1)}$. To see this, let $(a^i_n)$ be the matrix associated to $p$. Note that as $p\notin g(\mathcal{S}^{(2N+1)})$ we have $a_0^0=0$ and $a_1^1+a_1^{-1}\le 1$. Suppose for example that $a_1^1=1$. Then we have $p\in f_q(\mathcal{S}^{(2N+1)})$ for $q={\scriptstyle (0,\ldots,0,1,0,\ldots,0)}$. We assume from now that $a^0_0=a^1_1=a^{-1}_1=0$.
Equivalently, the stopping time $\tau$ described by $(a_n^i)$ satisfies $\p(\tau\ge 2)= 1$, which we can yet rephrase to say that two paths arrive in zero after two steps.
We now \emph{try} and construct a matrix $(\td{a}^i_n)$ to correspond to an embedding of $4p$ --- although this will not be strictly possible, it will determine the value of $q$ we will need so that $4(p-q)\in \mathcal{S}^{(2N+1)}$.
More precisely, define $\td{k}_0^i=0$ for all $i$, $\td{k}_1^i=0$ for all $i\neq 0$, $\td{k}_1^0=1$, and let $\td{a}^i_n=\max\{a^i_n,\td{k}^i_n\}$ where
\begin{equation}
\label{eqn:tdksystem}
\left\{
\begin{array}{lcl}
    \td{k}^{2i+1}_{n+1}&=& \td{k}^{2i}_{n}-\td{a}^{2i}_{n}+\td{k}^{2(i+1)}_{n}-\td{a}^{2(i+1)}_{n},\ n\ge 1,\ 2i+1\in [-N,N],\\
    \td{k}^{2i}_{n+1}&=&
\td{k}^{2i+1}_{n+1}-\td{a}^{2i+1}_{n+1}+\td{k}^{2i-1}_{n+1}-\td{a}^{2i-1}_{n+1},\ n\ge 1,\ 2i\in [-N,N].\\
\end{array}%
\right.
\end{equation}
Put $\td{p}^i=2^{(i\bmod 2)}\sum_{j=0}^\infty 4^{-j}\td{a}^i_j$ and $q=p-\td{p}$. From the construction, both $\td{p}$ and $q$ are elements of $\mathcal{S}^{(2N+1)}$ since their associated matrices are respectively $(\td{a}^i_n)$ and $(a^i_n-\td{a}^i_n)$. It is also clear that $q\in \mathcal{W}$ since we put explicitly $\td{k}_1^0=1$ as if `something else' stopped one of the two paths reaching zero after two steps.\\
We conclude that
\begin{eqnarray*}
\mathcal{S}^{(2N+1)}&=&\bigcup_{q\in \mathcal{W}}(\frac{1}{4}\mathcal{S}^{(2N+1)}+q)\cup g(\mathcal{S}^{(2N+1)})\\
   &=& \left(\frac{1}{4}\mathcal{S}^{(2N+1)}+\mathcal{W}\right)\cup g(\mathcal{S}^{(2N+1)})=f\left(\mathcal{S}^{(2N+1)}\right). \end{eqnarray*}
We would like to conclude that $f$ is a contraction and $\mathcal{S}^{(2N+1)}$ is its unique fixed point. To this end we need to show that $\mathcal{S}^{(2N+1)}$ and $\mathcal{W}$ are closed and thus compact (since both are bounded). Indeed, as Minkowski's sum of two compact sets is again compact, the mapping $f$ defined via $f(A)=(A/4+\mathcal{W})\cup g(A)$ is then a contraction on closed subsets of $[0,1]^{2N+1}$ and $\mathcal{S}^{(2N+1)}$ is its unique fixed point.\\
We show first that $\mathcal{S}^{(2N+1)}$ is closed. Consider a sequence $p_j\to p$, as $j\to\infty$, with $p_j\in \mathcal{S}^{(2N+1)}$. With each $p_j$ we have the associated matrix $(a^i_n(p_j))$, $|i|\le N$, $n\ge 0$. For a point $q\in \mathcal{S}^{(2N+1)}$ and its associated matrix $(a_n^i(q))$ we have (by the Optional Stopping Theorem) $\sum_{n=0}^\infty a^i_n(j)4^{-n}\le \frac{N+1}{N+1+|i|}$. In consequence, for any fixed depth $m\ge 1$, the set of matrices $\{(a_n^i(q)): |i|\le N, n\le m, q\in \mathcal{S}^{(2N+1)}\}$ is finite. We can therefore choose a subsequence $p_{g_j}\to p$ with the same matrix representation up to the depth $m$:  \begin{equation}\label{eqn:gettingmatrix}
a^i_n(p_{g_j})=a^i_n(p_{g_l}),\quad j,l\ge 0,\quad  n\le m.
\end{equation}
 We can then iterate the procedure. We can choose again a subsequence of the sequence $p_{g_j}$, such that \eqref{eqn:gettingmatrix} is verified for all $n\le 2 m$, then for $n\le 4 m$ and so on. In this way we obtain a matrix $A=(a_n^i)$ and a sequence $q_j\to p$ such that for any $d\ge 1$, $a_n^i=a_n^i(q_j)$ for all $|k|\le N$, $n\le d$ and $j\ge d$. In particular, the matrix $A$ satisfies $a_n^i\le k_n^i$ with $(k^i_n)$ defined via \eqref{eqn:ksystem}. Furthermore, we have
\begin{eqnarray}
    \sum_{n=0}^\infty 4^{-n}a^{i}_n&=&\lim_{d\to\infty}
\sum_{n=0}^d 4^{-n}a^{i}_n=\lim_{d\to\infty}\sum_{n=0}^d 4^{-n}a^{i}_n(q_d)\nonumber\\
&=&\lim_{d\to\infty}\left( 2^{-i\bmod 2}q^i_d-\sum_{n=d+1}^\infty 4^{-n}a^{i}_n(q_d)\right)=
2^{-i\bmod 2}p^i\ .\label{eq:gettingatoms}
\end{eqnarray}
To justify the last equality first note that $q_d\to p$ and so $q^i_d\to p^i$ as $d\to\infty$. Secondly, define $H_{N}=\inf\{n: X_n\notin [-N,N]\}$ and observe the upper bound $\sum_{n=d}^\infty 4^{-n}a^{i}_n(q_d)\le \p(H_N\ge d)\to 0$, as $d\to \infty$, since $\e H_N= (N+1)^2<\infty$.\\
Finally, $\mathcal{W}$ is clearly closed by its definition and the fact that $\mathcal{S}^{(2N+1)}$ is closed.
\end{proof}
To understand entirely the set $\mathcal{M}_0^{UI}$ it rests to describe its elements with unbounded support. To this end consider first $\mu\in\mathcal{M}$ any probability measure on $\zr$. Theorem \ref{thm:existence_min} implies existence of a minimal stopping time $\tau$ such that $X_\tau\sim\mu$. Let $\tau_N=\tau\land H_N$. Naturally $\tau_N\to \tau$ as $N\to \infty$ and thus $X_{\tau_N}\to X_\tau$ a.s.. Furthermore, as $(X_{\tau_N\land n}:n\ge 0)$ is a UI martingale, the measure $\mu_N$, the law of $X_{\tau_N}$, is an element of $\mathcal{S}^{(2N+1)}$. Thus if we consider the set of all measures with bounded support which can be embedded via UI stopping times
\begin{equation}\label{eqn:Sclosure}
\mathcal{S}^\infty=\bigcup_{N\ge 1}\mathcal{S}^{(2N+1)}\subset \mathcal{M}_0^{UI}\quad\textrm{then}\quad\overline{\mathcal{S}^\infty}=\mathcal{M},
\end{equation}
where the closure is taken in the topology of weak convergence. \\
In order to study closures in different topologies we identify for the rest of this paragraph, sets of measures with sets of random variables, so that $\mathcal{S}^\infty=\{X_\tau: \exists N\, \tau\le H_N\}$, with $\tau$ a stopping time, and likewise for $\mathcal{M}_0^{UI}$, $\mathcal{M}_0$ and $\mathcal{M}$. Furthermore, introduce the $L^p$ subsets of the set $\mathcal{M}^{UI}_0$:
\begin{equation*}\label{eqn:LpUIset}
    \mathcal{M}^{UI,p}_0=\Big\{X\in \mathcal{M}_0^{UI}: \e |X|^p<\infty\Big\},\quad p\ge 1.
\end{equation*}
 Then the following proposition holds.
\begin{proposition}
For any $p\ge 1$, $\mathcal{M}^{UI,p}_0$ is the closure of $\mathcal{S}^\infty$ in the $L^p$ norm:
\begin{equation}\label{eq:L2closure}
 \overline{\mathcal{S}^\infty}^{L^p}= \mathcal{M}^{UI,p}_0.
\end{equation}
\end{proposition}
\begin{proof}
We prove first the inclusion "$\subset$".
Suppose that a sequence $X_{\rho_N}$ in $\mathcal{S}^\infty$ converges in $L^p$, $p\ge 1$, to some variable $X$. We can then replace $\rho_N$ with $\tau_N=\min\{\rho_K: K\ge N\}$ which is an increasing sequence of stopping times, which thus converges to a stopping time: $\tau_N\nearrow \tau$ \as{}. Further, since
\[
|X_n| - \sum_{k=0}^{n-1} \mathbf{1}_{\{X_k=0\}}
\]
is a martingale, we have
\[
\e |X_{\tau_N}| = \e\left( \sum_{k=0}^{\tau_N-1} \mathbf{1}_{\{X_k=0\}}\right).
\]
Noting that the left hand side is bounded since $\tau_N \le \rho_N$ and therefore $\e |X_{\tau_N}| \le \e |X_{\rho_N}|$, we obtain
\[
\e\left( \sum_{k=0}^{\tau} \mathbf{1}_{\{X_k=0\}}\right) < \infty,
\]
and from the recurrence of the random walk we can deduce that $\tau < \infty$ \as{}.\footnote{This is an analogue of an argument used originally in \cite{cox_hobson2} in the continuous setting.} In particular, we can now make sense of $X_\tau$. Therefore $X_{\tau_N}\to X_\tau$ \as{} and in $L^p$ as $N\to \infty$, and so \emph{a fortiori} $X_\tau=X$ a.s. In consequence,  $(X_{\tau_N}: N\ge 1)$ is a uniformly integrable martingale. Furthermore, for every $N\ge 1$, $(X_{\tau_N\land t}: t\ge 0)$ is also a UI martingale. We have thus
\begin{equation*}
%\label{eq}
X_{t\land \tau_N}=\e \Big[X_{\tau_N}\Big| \mathcal{F}_{t\land \tau_N}\Big]=\e \Big[\e[X_\tau |\mathcal{F}_{\tau_N}]\Big| \mathcal{F}_{t\land \tau_N}\Big]=\e \Big[X_{\tau}\Big| \mathcal{F}_{t\land \tau_N}\Big]
\end{equation*}
and taking the limit as $N\to \infty$ we see that $X_{t\land \tau}=\e[X_\tau|\mathcal{F}_{t\land\tau}]$ a.s. (note that $\e |X_\tau|\le\infty$). This proves that $X=X_\tau\in \mathcal{M}^{UI,p}_0$.\\
The converse is easier. Let $X_\tau\in\mathcal{M}^{UI,p}_0$ and put $\tau_N=\tau\land H_N$. Then $X_{\tau_N}=\e[X_\tau|\mathcal{F}_{\tau_N}]$ converges a.s.\ and in $L^1$ to $X_\tau$ as $N\to \infty$. The convergence actually holds in $L^p$ as $\sup_N \e |X_{\tau_N}|^p=\e |X_\tau|^p<\infty$ (cf.\ Revuz and Yor \cite[Thm II.3.1]{MR2000h:60050}). Naturally, $X_{\tau_N}\in S^\infty$ and thus $X_\tau \in \overline{\mathcal{S}^\infty}^{L^p}$.
\end{proof}

\section{Conclusions and Further problems}
We have studied the Skorokhod embedding problem for the simple symmetric random walk and the relations between various classes of stopping times. In particular, we have seen that --- unlike in the Brownian motion setup --- the classes of uniformly integrable and minimal stopping times are not equal. The latter allows us to construct an embedding for any centered target measure. %The former restricts the class of admissible measures, and further it possesses a complex fractal structure of the set of UI embeddable measures with bounded support.\\
The former restricts the class of admissible measures, and in fact we show that the set of UI embeddable measures with bounded support has a complex fractal structure.\\
Our study answered thus the questions we have asked ourselves at the beginning. We would like to stress however, that it raised at least as many new questions, which seem interesting to us. We review few of them to end the paper.\\
It would be interesting to calculate the Lebesgue measure of $\mathcal{S}^{(2N+1)}$ and to study further its structure; we have not considered the local dimension of the sets in higher dimensions. We would also like to understand the relationship (\eg{} as projections) between the sets for different values of $N$.\\
%We would also like to see a nice characterization of the measures in $\mathcal{M}^{UI}_0$ with unbounded support.\\
As far as the Skorokhod embedding is concerned, we note that we have not given an explicit construction for every $\mu \in \mathcal{S}^\infty$ or $\mu \in \mathcal{M}_0^{UI}$. It seems a hard but interesting goal. We have not really uncovered the meaning of minimality of stopping times for the random walk. We show that it is very different from the continuous martingale setup but we have not devised any criterion, given in terms of the stopped process, to decide whether a given stopping time is minimal.
Understanding minimality of stopping times and extending the results to arbitrary discontinuous martingales remains an open problem.

\bibliographystyle{abbrv}
%\bibliography{../../bib/bibliografia}
\bibliography{../../bibliografia}
\end{document}